\documentclass[preprint,11pt]{elsarticle}

\usepackage{amssymb}
\usepackage{amsmath}
\usepackage{amsfonts}
\usepackage{amssymb}
\usepackage{indentfirst,latexsym,bm}
\usepackage{amsthm}
\usepackage{color,xcolor}
\usepackage[all]{xy}
\input{mathrsfs.sty}
\newtheorem{theorem}{Theorem}[section]
\newtheorem{lemma}[theorem]{Lemma}
\newtheorem{corollary}[theorem]{Corollary}

\theoremstyle{definition}
\newtheorem{definition}{Definition}[section]

\theoremstyle{definition}
\newtheorem{example}{Example}[section]

\theoremstyle{remark}
\newtheorem{remark}{Remark}[section]
\theoremstyle{question}

\numberwithin{equation}{section}

\journal{XXX}

\begin{document}

\begin{frontmatter}



\title{$C^*$-isomorphisms associated with two projections on a Hilbert $C^*$-module}
\author[rvt]{Chunhong Fu}
\ead{fchlixue@163.com}
\author[rvt]{Qingxiang Xu\corref{cor1}}
\ead{qingxiang$\_$xu@126.com}
\cortext[cor1]{Corresponding author}
\author[rvt]{Guanjie Yan}
\ead{m15578192102@163.com}
\address[rvt]{Department of Mathematics, Shanghai Normal University, Shanghai 200234, PR China}

\begin{abstract}
 Motivated by two norm equations used to characterize the Fried\-richs angle, this paper studies $C^*$-isomorphisms associated with two  projections by introducing the matched triple and the semi-harmonious pair of projections.
A triple $(P,Q,H)$ is said to be matched if $H$ is a Hilbert $C^*$-module, $P$ and $Q$ are projections on $H$  such that
 their infimum $P\wedge Q$ exists as an element of $\mathcal{L}(H)$, where $\mathcal{L}(H)$ denotes the set of all  adjointable operators on $H$. The $C^*$-subalgebras of $\mathcal{L}(H)$ generated by elements in $\{P-P\wedge Q, Q-P\wedge Q, I\}$ and
 $\{P,Q,P\wedge Q,I\}$ are denoted by $i(P,Q,H)$ and $o(P,Q,H)$, respectively.
It is proved that each faithful representation $(\pi, X)$ of $o(P,Q,H)$ can induce a faithful representation $(\widetilde{\pi}, X)$ of $i(P,Q,H)$
 such that
 \begin{align*}&\widetilde{\pi}(P-P\wedge Q)=\pi(P)-\pi(P)\wedge \pi(Q),\\
&\widetilde{\pi}(Q-P\wedge Q)=\pi(Q)-\pi(P)\wedge \pi(Q).
\end{align*}
 When $(P,Q)$ is semi-harmonious, that is,  $\overline{\mathcal{R}(P+Q)}$ and $\overline{\mathcal{R}(2I-P-Q)}$ are both orthogonally complemented in $H$, it is shown that  $i(P,Q,H)$ and $i(I-Q,I-P,H)$ are unitarily equivalent via a unitary operator in $\mathcal{L}(H)$. A counterexample is constructed, which shows that the same may be not true when $(P,Q)$ fails to be semi-harmonious. Likewise, a  counterexample is constructed such that
 $(P,Q)$ is semi-harmonious, whereas $(P,I-Q)$ is not semi-harmonious. Some additional examples indicating new phenomena of adjointable operators acting on Hilbert $C^*$-modules are also provided.

\end{abstract}


\end{frontmatter}



\section{Introduction} \label{sec:introduction}

Let $P$ and $Q$ be two projections on a Hilbert space $X$. Their infimum $P\wedge Q$ is the projection from $X$ onto $\mathcal{R}(P)\cap \mathcal{R}(Q)$, which can be obtained by taking the limit of  $\{(PQP)^n\}_{n=1}^\infty$ in the strong operator topology  \cite[Lemma~22]{von Neumann}. The cosine of the Friedrichs angle \cite{Friedrichs}  between $M=\mathcal{R}(P)$ and $N=\mathcal{R}(Q)$ is denoted by $c(M,N)$, and  can be calculated  as
 \begin{equation*}c(M,N)=\left\|(P-P\wedge Q)(Q-P\wedge Q)\right\|.\end{equation*}
The characterization of $c(M,N)=c(N^\bot,M^\bot)$ given in \cite[Section~2]{Deutsch}  yields
 \begin{align}\label{equ:motivation equation-1}&\big\|(P-P\wedge Q)(Q-P\wedge Q)\big\|\ \nonumber\\
&\quad =\Big\| \big[I-Q-(I-Q)\wedge (I-P)\big]\big[I-P-(I-Q)\wedge (I-P)\big]\Big\|.
\end{align}

The Friedrichs angle has also been studied in the setting of $C^*$-algebras. To deal with the Friedrichs angle associated with two projections $P$ and $Q$ in a $C^*$-algebra $\mathfrak{A}$, one approach employed in \cite{Anoussis} is
to embed $\mathfrak{A}$ into its enveloping von Neumann algebra $\mathfrak{A}^{\prime\prime}$ via the universal representation $(\pi_u,H_u)$ of $\mathfrak{A}$,
and then to use the universal property of $\mathfrak{A}^{\prime\prime}$ \cite[Theorem~3.7.7]{Pedersen}. By identifying $\mathfrak{A}$ with $\pi_u(\mathfrak{A})$,  $P\wedge Q$ can be obtained in $\mathfrak{A}^{\prime\prime}$,
 and it is proved in \cite[Proposition~2.4]{Anoussis}
 that for every faithful representation $(\pi,X)$ of $\mathfrak{A}$,
 \begin{align}\label{equ:motivation equation-2}&\big\|(P-P\wedge Q)(Q-P\wedge Q)\big\|\nonumber\\
&\quad =\Big\|\big[\pi(P)-\pi(P)\wedge \pi(Q)\big]\big[\pi(Q)-\pi(P)\wedge \pi(Q)\big]\Big\|.
 \end{align}

 Hilbert $C^*$-modules are natural generalizations of Hilbert spaces, and every $C^*$-algebra can be regarded as a Hilbert $C^*$-module over itself in a natural way. The purpose of this paper is, in the framework of Hilbert $C^*$-modules, to give a deeper understanding of
 \eqref{equ:motivation equation-1} and \eqref{equ:motivation equation-2} via algebraic systems rather than on  products  of finitely many operators.

It is notable that a closed submodule of a Hilbert $C^*$-module may fail to be orthogonally complemented. In this paper much attention has been paid on this aspect. Let $\mathcal{L}(H)$ be the set of all adjointable operators on a Hilbert $C^*$-module $H$. For two projections $P,Q\in\mathcal{L}(H)$,  let  \begin{equation}\label{defn of R and N}\mathcal{R}=\mathcal{R}(Q)\cap \mathcal{R}(P)\quad\mbox{and}\quad  \mathcal{N}=\mathcal{N}(Q)\cap \mathcal{N}(P).
 \end{equation}
To check the validity of the Halmos' two projections theorem in the Hilbert $C^*$-module case, the term of the harmonious pair of projections is introduced in \cite[Section~4]{Luo-Moslehian-Xu}, and it is  shown later in \cite[Theorem~3.3]{Xu-Yan} that for every  pair $(P,Q)$ of projections, the Halmos' two projections theorem is valid  if and only if $(P,Q)$ is harmonious.
In view of the conditions stated in \cite[Lemma~5.4]{Luo-Moslehian-Xu} and \cite[Theorem~3.3]{Xu-Yan},  we make a definition as follows.

\begin{definition}\quad\label{defn of semi} A pair $(P,Q)$ of projections on a Hilbert $C^*$-module $H$ is said to be semi-harmonious if both $\overline{\mathcal{R}(P+Q)}$ and $\overline{\mathcal{R}(2I-P-Q)}$ are  orthogonally complemented in $H$. If $(P,Q)$ and $(P,I-Q)$ are
both semi-harmonious, then $(P,Q)$ is said to be harmonious.
\end{definition}
Let $\mathcal{R}$ and $\mathcal{N}$ be defined by \eqref{defn of R and N} for projections $P$ and $Q$ on a Hilbert $C^*$-module $H$.  Since $\overline{\mathcal{R}(2I-P-Q)}^\bot=\mathcal{R}$ and  $\overline{\mathcal{R}(P+Q)}^\bot=\mathcal{N}$,  a condition weaker than the  semi-harmony of $(P,Q)$ turns out to be the orthogonal complementarity of $\mathcal{R}$ and $\mathcal{N}$, which is necessary and sufficient to  make use of the notations $P_{\mathcal{R}}$ and $P_{\mathcal{N}}$ (the projections from $H$ onto $\mathcal{R}$ and $\mathcal{N}$, respectively). The example constructed in \cite[Section~3]{Luo-Moslehian-Xu} (see also the proof of Theorem~\ref{thm:semi-harmonious not implication}) shows that there exist projections $P$ and $Q$ on certain Hilbert $C^*$-module such that $\mathcal{R}=\mathcal{N}=\{0\}$, whereas $(P,Q)$ is not semi-harmonious. So, generally the meaningfulness of $P_{\mathcal{R}}$ and $P_{\mathcal{N}}$ does not imply the semi-harmony of $(P,Q)$.

Next, we introduce the matched triple as follows.

\begin{definition}\quad\label{defn of matched triple} A triple $(P,Q,H)$ is said to be matched if $H$ is a Hilbert $C^*$-module, $P$ and $Q$ are projections on $H$  such that
 $\mathcal{R}$ defined by \eqref{defn of R and N} is orthogonally complemented in $H$. In this case, the projection $P_{\mathcal{R}}$ is denoted by $P\wedge Q$. The $C^*$-subalgebras of $\mathcal{L}(H)$ generated by elements in $\{P-P\wedge Q, Q-P\wedge Q, I\}$ and
 $\{P,Q,P\wedge Q,I\}$ are denoted by $i(P,Q,H)$ and $o(P,Q,H)$,  and are called the inner algebra and the outer algebra, respectively.
\end{definition}

\begin{definition}\quad Two matched triple $(P_i,Q_i,H_i)$ $(i=1,2)$ are said to be innerly unitarily equivalent
if there exists a unitary operator $U: H_1\to H_2$ such that
\begin{align*}&U(P_1-P_1\wedge Q_1)U^*=P_2-P_2\wedge Q_2,\\
&U(Q_1-P_1\wedge Q_1)U^*=Q_2-P_2\wedge Q_2.
\end{align*}
\end{definition}

Recall that a pair $(\pi,X)$ is said to be a representation of a $C^*$-algebra $\mathfrak{A}$ if $X$ is a Hilbert space and $\pi: \mathfrak{A}\to \mathbb{B}(X)$ is a $C^*$-morphism, where $\mathbb{B}(X)$ denotes the set of all bounded linear operators on $X$.
\begin{definition}\quad\label{shd defn}Let $(P,Q,H)$ be a matched triple. A representation $(\pi,X)$ of $o(P,Q,H)$ is called an outer representation of $(P,Q,H)$. If furthermore a $C^*$-morphism $\widetilde{\pi}: i(P,Q,H)\to \mathbb{B}(X)$ can be induced such that $\widetilde{\pi}(I)=\pi(I)$, and
\begin{align*}&\widetilde{\pi}(P-P\wedge Q)=\pi(P)-\pi(P)\wedge \pi(Q),\\
&\widetilde{\pi}(Q-P\wedge Q)=\pi(Q)-\pi(P)\wedge \pi(Q),
\end{align*}
then $(\pi,X)$ is called an inner-outer representation of $(P,Q,H)$. When both $\pi$ and $\widetilde{\pi}$ are faithful, $(\pi,X)$ is called a faithful inner-outer representation of $(P,Q,H)$.
\end{definition}

\begin{remark}\quad Let  $(\pi, X)$ be an inner-outer representation of $(P,Q,H)$. It is notable that generally $\pi(P)\wedge \pi(Q)$ is taken in the von Neumann algebra $\Big[\pi\big[o(P,Q,H)\big]\Big]^{\prime\prime}$ rather than in the $C^*$-algebra $\pi\big[o(P,Q,H)\big]$, so it may happen that
$\pi(P\wedge Q)\ne \pi(P)\wedge \pi(Q)$.
\end{remark}

With the terms given as above, we list the main results of this paper  as follows:

1. There exist projections $P$ and $Q$ on certain  Hilbert $C^*$-module $H$ such that
$(P,Q)$ is semi-harmonious, whereas it fails to be harmonious (see Theorem~\ref{thm:semi-harmonious not implication}).

2.  For every semi-harmonious pair $(P,Q)$ of projections on a Hilbert $C^*$-module $H$, the matched triples $(P,Q,H)$ and $(I-Q,I-P,H)$ are innerly unitarily equivalent (see Theorem~\ref{technical thm}).

3.  There exist projections $P$ and $Q$ on certain  Hilbert $C^*$-module $H$ such that
the triples $(P,Q,H)$ and $(I-Q,I-P,H)$ are both matched, whereas they are not innerly unitarily equivalent (see Theorem~\ref{thm:unitary equivalence stops}).

4.  Every faithful outer representation of a matched triple is a faithful inner-outer representation (see Theorem~\ref{thm:algebras isomorphic-01}).

An application of Theorems~\ref{technical thm} and \ref{thm:algebras isomorphic-01}  will be illustrated in Corollary~\ref{cor:be applied}. Another application, as has been mentioned earlier, concerns a new insight into \eqref{equ:motivation equation-1} and \eqref{equ:motivation equation-2}. Let $P$ and $Q$ be two projections on a Hilbert $C^*$-module $H$ such that
$\mathcal{R}$ and $\mathcal{N}$ defined by \eqref{defn of R and N} are orthogonally complemented in $H$.
By Theorem~\ref{thm:algebras isomorphic-01}, we will see that
 each faithful unital representation $(\pi,X)$  of $\mathcal{L}(H)$ can induce unital $C^*$-isomorphisms $\widetilde{\pi}_1: i(P,Q,H)\to i\big(\pi(P),\pi(Q),X\big)$ and $\widetilde{\pi}_2: i(I-Q,I-P,H)\to i\big(I-\pi(Q),I-\pi(P),X\big)$
  such that
  \begin{align*}&\widetilde{\pi}_1(P-P\wedge Q)=\pi(P)-\pi(P)\wedge \pi(Q),\\
  &\widetilde{\pi}_1(Q-P\wedge Q)=\pi(Q)-\pi(P)\wedge \pi(Q),\\
  &\widetilde{\pi}_2\big[I-Q-(I-Q)\wedge (I-P)\big]=I-\pi(Q)-\pi(I-Q)\wedge \pi(I-P),\\
  &\widetilde{\pi}_2\big[I-P-(I-Q)\wedge (I-P)\big]=I-\pi(P)-\pi(I-Q)\wedge \pi(I-P).
\end{align*}
Thus, $\|\widetilde{\pi}_1(x)\|=\|x\|$ for every $x\in i(P,Q,H)$. Specifically, if we put $x=(P-P\wedge Q)(Q-P\wedge Q)$, then \eqref{equ:motivation equation-2} is obtained.

Note that $\big(\pi(P), \pi(Q)\big)$ is a pair of projections acting  on a Hilbert space, so it is harmonious. Hence, by Theorem~\ref{technical thm} there exists a unitary operator $U\in\mathbb{B}(X)$ such that
\begin{align*}&U\big[\pi(P)-\pi(P)\wedge \pi(Q)\big]U^*=I-\pi(Q)-\pi(I-Q)\wedge \pi(I-P),\\
&U\big[\pi(Q)-\pi(P)\wedge \pi(Q)\big]U^*=I-\pi(P)-\pi(I-Q)\wedge \pi(I-P).
\end{align*}
Thus, a $C^*$-isomorphism $\rho: i(P,Q,H)\to  i(I-Q,I-P,H)$ can be constructed as
\begin{equation*}\label{defn of rho-01}\rho(x)=\left(\widetilde{\pi}_2\right)^{-1}U\widetilde{\pi}_1(x)U^*, \quad\forall x\in i(P,Q,H).
\end{equation*}
Therefore, $\|\rho(x)\|=\|x\|$ for every  $x\in i(P,Q,H)$. Likewise, if we take $x=(P-P\wedge Q)(Q-P\wedge Q)$, then \eqref{equ:motivation equation-1} is obtained. So,
a substantive generalization of \eqref{equ:motivation equation-1}  and \eqref{equ:motivation equation-2} has been made.

The paper is organized as follows. The main purpose of Section~\ref{sec:Semi-harmonious} is to construct two projections $P$ and $Q$ such that $(P,Q)$ is semi-harmonious, whereas
 it fails to be harmonious. Section~\ref{sec:unitary equivalence} focuses on the construction of the unitary operator $U$ satisfying \eqref{equation of u equivalence-1} and \eqref{equation of u equivalence-2}. Section~\ref{sec:norm equations}
is devoted to the study of the faithful inner-outer representation of a matched triple.

\section{Semi-harmonious pairs of projections}\label{sec:Semi-harmonious}

Throughout the rest of this paper,  $\mathbb{N}$, $\mathbb{Z}_+$ and $\mathbb{C}$ are the sets of all positive integers, non-negative integers and complex numbers, respectively. Unless otherwise specified, $\mathfrak{A}$ is a $C^*$-algebra,  $E, H$ and $K$ are Hilbert $\mathfrak{A}$-modules \cite{Lance,Paschke}.
The set of  all adjointable operators from $H$ to $K$ is denoted by $\mathcal{L}(H,K)$. Given $A\in\mathcal{L}(H,K)$, the adjoint operator, the range and the null space of $A$ are denoted by $A^*$,
$\mathcal{R}(A)$ and $\mathcal{N}(A)$, respectively. Let $|A|$ designate the square root of $A^*A$.
In case $H=K$, $\mathcal{L}(H,K)$ is abbreviated to $\mathcal{L}(H)$, whose subset consisting of all positive elements is denoted by $\mathcal{L}(H)_+$. The unit of $\mathcal{L}(H)$ (namely, the identity operator on $H$) is denoted by $I_H$, or simply by $I$ when no ambiguity arises.
An operator $P\in \mathcal{L}(H)$
is said to be a projection if $P = P^*= P^2$.  The set of all projections on  $H$ is denoted by $\mathcal{P}(H)$.

Let $M$ be a closed submodule of $H$. Clearly, there exists at most a projection in $\mathcal{L}(H)$, written $P_M$, such that $\mathcal{R}(P_M)=M$. It can be easily verified that $P_M$ exists if and only if
$H=M+M^\bot$, where $$M^\bot=\left\{x\in H: \langle x,y\rangle=0,\forall y\in M\right\}.$$ In this case,
$M$ is said to be orthogonally complemented in $H$.

To construct semi-harmonious pairs of projections, we need a couple of lemmas.
\begin{lemma}\quad\label{lem:Range Closure of T alpha and T}{\rm \cite[Proposition~2.9]{Liu-Luo-Xu}} For every $T\in \mathcal{L}(H)_+$ and $\alpha>0$, we have  $\overline{\mathcal{R}(T)}=\overline{\mathcal{R}(T^{\alpha})}$.
\end{lemma}

\begin{lemma}\quad\label{lem:Range closure of TT and T} {\rm\cite[Proposition 3.7]{Lance}} For every $T\in\mathcal{L}(H,K)$, we have $\overline{\mathcal{R}(T)}=\overline{\mathcal{R}(TT^*)}$.
\end{lemma}

\begin{lemma}\quad\label{lem:rang characterization-1}{\rm \cite[Proposition~2.7]{Liu-Luo-Xu}} Let $B,C\in\mathcal{L}(E,H)$ be such that $\overline{\mathcal{R}(B)}=\overline{\mathcal{R}(C)}$. Then for every $A\in\mathcal{L}(H,K)$, we have $\overline{\mathcal{R}(AB)}=\overline{\mathcal{R}(AC)}$.
\end{lemma}

An approach to construct  semi-harmonious pairs of projections reads as follows.

\begin{theorem}\quad  \label{thm:invariant submodule M}For every $P,Q\in \mathcal{P}(H)$, let $M\subseteq H$ be defined by
\begin{equation}\label{defn of M}M=\overline{\mathcal{R}\big[(P+Q)(2I-P-Q)\big]}.\end{equation}
Then  $P|_M$ and $Q|_M$ are projections on $M$ such that $(P|_M,Q|_M)$ is semi-harmonious, where $P|_M$ and $Q|_M$ are the restrictions of $P$ and $Q$ on $M$, respectively.
\end{theorem}
\begin{proof}[\bf Proof]To simplify the notation, we put \begin{equation}\label{defn of A and G}A=P+Q \quad\mbox{and}\quad  G=\mathcal{R}\big[A(2I-A)\big].\end{equation}
Since $G$ defined as above is the range of an adjointable operator, its closure $M$ is a Hilbert $\mathfrak{A}$-module.

For every $T\in\mathcal{L}(H)$ and $X\subseteq H$, let $TX=\big\{Tx: x\in X\big\}$. Then
$TX=\mathcal{R}(T|_X)$, and the boundedness of $T$  gives $\overline{TX}=\overline{T\overline{X}}$. Hence
\begin{equation}\label{trivial equation of closures}\overline{AM}=\overline{AG}\quad\mbox{and}\quad \overline{(2I-A)M}=\overline{(2I-A)G}.
\end{equation}
Since  $A$ and $2I-A$ are positive and commutative, we may combine \eqref{defn of A and G} and \eqref{trivial equation of closures} with
Lemmas~\ref{lem:Range Closure of T alpha and T} and \ref{lem:rang characterization-1} to get
\begin{align}\label{overline is full-01}\overline{\mathcal{R}(P|_M+Q|_M)}=&\overline{AM}=\overline{AG}=\overline{\mathcal{R}\big[(2I-A)A^2\big]}\nonumber\\=&\overline{\mathcal{R}\big[(2I-A)A\big]}=M.
\end{align}
Similarly,
\begin{align}\label{overline is full-02}\overline{\mathcal{R}(2I_M-P|_M-Q|_M)}=&\overline{(2I-A)M}=\overline{(2I-A)G}=\overline{\mathcal{R}\big[A(2I-A)^2\big]}\nonumber\\
=&\overline{\mathcal{R}\big[A(2I-A)\big]}=M.
\end{align}
Furthermore, direct computations yield
\begin{align*}
&P(2I-A)A=P(I-Q)P=A(2I-A)P,\\
&Q(2I-A)A=Q(I-P)Q=A(2I-A)Q,
\end{align*}
which lead clearly to $PM\subseteq M$ and $QM\subseteq M$. Consequently, $P|_M$ and $Q|_M$ are  projections on $M$. In view of \eqref{overline is full-01} and \eqref{overline is full-02}, we conclude that  $(P|_M,Q|_M)$ is semi-harmonious.
\end{proof}

\begin{theorem}\quad  \label{thm:semi-harmonious not implication}There exist projections $P$ and $Q$ on certain Hilbert $C^*$-module such that $(P,Q)$ is semi-harmonious, whereas it fails to be harmonious.
\end{theorem}
\begin{proof}[\bf Proof]We follow the line initiated in \cite[Section~3]{Manuilov-Moslehian-Xu} and modified in \cite[Section~3]{Luo-Moslehian-Xu}.
Let $M_2(\mathbb{C})$ and $I_2$ be the set of all $2\times 2$ complex matrices and the identity matrix in $M_2(\mathbb{C})$, respectively.
Denote by $\|\cdot\|$ the operator norm on $M_2(\mathbb{C})$.
Let $\mathfrak{A}=C\big([0,1];M_2(\mathbb{C})\big)$ be the set of all continuous matrix-valued functions from $[0,1]$ to $M_2(\mathbb{C})$. For $x\in \mathfrak{A}$ and $t\in [0,1]$, we put
\begin{align*}x^*(t)=\big(x(t)\big)^*\quad \mbox{and}\quad \|x\|=\max_{0\le s\le 1}\|x(s)\|.\end{align*}
With the $*$-operation as above and the usual algebraic operations, $\mathfrak{A}$ is a unital $C^*$-algebra, which is also a Hilbert $\mathfrak{A}$-module with the inner-product given by
$$\langle x,y\rangle=x^*y,\quad \mbox{for $x,y\in \mathfrak{A}$}.$$
Let $e$ be the unit of $\mathfrak{A}$, that is, $e(t)=I_2$ for every $t\in [0,1]$.
It is known  that $\mathfrak{A}\cong \mathcal{L}(\mathfrak{A})$ via $a\to L_a$ \cite[Section~3]{Luo-Moslehian-Xu}, where
$L_a(x)=ax$ for $a,x\in \mathfrak{A}$. For simplicity, we identify $\mathcal{L}(\mathfrak{A})$ with $\mathfrak{A}$ and set
\begin{equation}\label{defn of s t and c t}c_t=\cos\frac{\pi}{2}t\quad \mbox{and}\quad s_t=\sin\frac{\pi}{2}t,\quad \mbox{for $t\in[0,1]$}.\end{equation} Let $P,Q\in \mathfrak{A}$ be projections determined by the matrix-valued functions
\begin{equation}\label{equ:defn of P t and Q t}P(t)\equiv\left(\begin{matrix}1&0\\0&0\end{matrix}\right)\quad \mbox{and}\quad Q(t)=\left(\begin{matrix}c_t^2&s_tc_t\\s_tc_t&s_t^2\end{matrix}\right), \quad \mbox{for $t\in [0,1]$}.\end{equation}
It is shown in  \cite[Section~3]{Luo-Moslehian-Xu} that
\begin{equation}\label{equ:intersections of null spaces are zeros}\mathcal{R}(P)\cap\mathcal{R}(Q)=\mathcal{R}(P)\cap\mathcal{N}(Q)=\mathcal{N}(P)\cap\mathcal{R}(Q)=\mathcal{N}(P)\cap\mathcal{N}(Q)=\{0\}.
\end{equation}

Now, let $H=\mathfrak{A}$ and let $M$ be defined by \eqref{defn of M}. According to Theorem~\ref{thm:invariant submodule M},
$(P|_M,Q|_M)$ is semi-harmonious. In what follows, we prove that $(P|_M,Q|_M)$ is not harmonious.

Direct computation yields
$$(P+Q)(2I-P-Q)=P+Q-PQ-QP=(P-Q)(P-Q)^*,$$
which leads by \eqref{defn of M} and Lemma~\ref{lem:Range closure of TT and T} to
\begin{equation}\label{equ:another from of M}M=\overline{\mathcal{R}(P-Q)}.\end{equation}
Utilizing \eqref{equ:defn of P t and Q t} we obtain
$P-Q=au$, where $a,u\in H$ are determined by
$$a(t)=\left(
         \begin{array}{cc}
           s_t & 0 \\
           0 & s_t \\
         \end{array}
       \right), \quad u(t)=\left(
       \begin{array}{cc}
         s_t & -c_t \\
         -c_t & -s_t \\
       \end{array}
     \right),\quad t\in [0,1].
$$
Since $u$ is a unitary and $M$ can be represented  by \eqref{equ:another from of M},  we have
$M=\overline{\mathcal{R}(a)}$. Specifically,
\begin{equation*}a=ae\in M,\quad \overline{\mathcal{R}(P|_M+I_{M}-Q|_M)}=\overline{\mathcal{R}(T)},\end{equation*} where
$T=(P+I-Q)a\in H$. For  any $x\in H$ with
$x(t)=\big(x_{ij}(t)\big)_{1\le i,j\le 2}$, it is easy to verify that
\begin{align*}(Tx)(1)=\left(
                                                 \begin{array}{cc}
                                                   2x_{11}(1) & 2x_{12}(1)\\
                                                   0 & 0 \\
                                                 \end{array}
                                               \right),
\end{align*}
hence
$$\|Tx-a\|\ge \|(Tx)(1)-a(1)\|=\left\|\left(
                                                 \begin{array}{cc}
                                                   2x_{11}(1)-1 & 2x_{12}(1)\\
                                                   0 & -1 \\
                                                 \end{array}
                                               \right)\right\|\ge 1,$$
which implies that $a\notin \overline{\mathcal{R}(T)}$. Furthermore, by \eqref{equ:intersections of null spaces are zeros} we have
\begin{align*}\overline{\mathcal{R}(P|_M+I_{M}-Q|_M)}^\bot=\mathcal{N}(P|_M)\cap \mathcal{R}(Q|_M)\subseteq \mathcal{N}(P)\cap \mathcal{R}(Q)=\{0\}.
\end{align*}
This shows $$a\notin \overline{\mathcal{R}(P|_M+I_{M}-Q|_M)}+\overline{\mathcal{R}(P|_M+I_{M}-Q|_M)}^\bot,$$ whereas $a\in M$.
So  $\overline{\mathcal{R}(P|_M+I_{M}-Q|_M)}$ is not orthogonally complemented in $M$.
\end{proof}

\begin{remark}\quad It is notable that there exist projections $P$ and $Q$ such that $\overline{\mathcal{R}(P+Q)}$ is orthogonally complemented, whereas $\overline{\mathcal{R}(2I-P-Q)}$ fails to be orthogonally complemented. We provide such an example as follows.
\end{remark}

\begin{example}\quad\label{ex:counter example-01}{\rm Let $\mathfrak{A}=H=C\big([0,1];M_2(\mathbb{C})\big)$ and let $P,Q\in\mathcal{P}(H)$ be as in   the proof of Theorem~\ref{thm:semi-harmonious not implication}.  Put $H_0=\overline{\mathcal{R}(P+Q)}$, $P_{0}=P|_{H_0}$ and $Q_{0}=Q|_{H_0}$.  From \cite[Lemma~2.3]{Luo-Moslehian-Xu} we have
$H_0=\overline{\mathcal{R}(P)+\mathcal{R}(Q)}$, which means that
$\mathcal{R}(P)\subseteq H_0$, hence $P_0 H_0\subseteq H_0$. Consequently, $P_0$ is a projection on $H_0$. Similarly, we have $Q_0\in\mathcal{P}(H_0)$. In view of \eqref{equ:intersections of null spaces are zeros},  we get
\begin{equation}\label{two intersections zero-1}\mathcal{R}(P_0)\cap \mathcal{R}(Q_0)=\mathcal{N}(P_0)\cap \mathcal{N}(Q_0)=\{0\}.\end{equation}
According to  Lemma~\ref{lem:Range Closure of T alpha and T},  we have
$$H_0=\overline{\mathcal{R}\big[(P+Q)^2\big]}\subseteq \overline{\mathcal{R}(P_0+Q_0)}\subseteq H_0.$$
As a result, we arrive at
$$H_0=\overline{\mathcal{R}(P_0+Q_0)}=\overline{\mathcal{R}(P_0+Q_0)}+\mathcal{N}(P_0)\cap\mathcal{N}(Q_0).$$
This shows the orthogonal complementarity of
$\overline{\mathcal{R}(P_{0}+Q_{0})}$ in $H_0$.

Let $F=(2I-P-Q)(P+Q)$. Clearly, $F=(I-P)Q+(I-Q)P$, so by \eqref{equ:defn of P t and Q t} we have
$$F(t)=\left(
         \begin{array}{cc}
           s_t^2 &  \\
            & s_t^2 \\
         \end{array}
       \right),\quad \forall t\in [0,1],$$
which implies that for every $x\in \mathfrak{A}$,
$(Fx)(0)=F(0)x(0)=\left(
                    \begin{array}{cc}
                      0 &  \\
                       & 0 \\
                    \end{array}
                  \right)$. Hence
$$\|P-Fx\|\geq \|P(0)-F(0)x(0)\|=\left\|\left(
                                \begin{array}{cc}
                                  1 & 0 \\
                                  0 & 0 \\
                                \end{array}
                              \right)
\right\|=1.$$
Due to the definition of $F$ and the observation of \eqref{two intersections zero-1}, we conclude that
$$P\notin\overline{\mathcal{R}(F)}=\overline{\mathcal{R}(2I_{H_0}-P_{0}-Q_{0})}+\mathcal{R}(P_0)\cap \mathcal{R}(Q_0).$$
On the other hand, $P=PI\in\mathcal{R}(P)\subseteq H_0$. Therefore, $\overline{\mathcal{R}(2I_{H_0}-P_{0}-Q_{0})}$ is not orthogonally complemented in $H_0$.
}\end{example}

\section{Unitary equivalences associated with two projections}\label{sec:unitary equivalence}\label{sec:unitary equivalence}

In this section, we deal with  unitary equivalences associated with two projections. We begin with a known result as follows.

\begin{lemma}\quad\label{charac of orth-PQ}{\rm \cite[Lemma~4.1]{Luo-Moslehian-Xu}}\  Let $P,Q\in\mathcal{P}(H)$ be such that $\overline{\mathcal{R}(I-Q+P)}$ is
orthogonally complemented in $H$. Then $\overline{\mathcal{R}(QP)}$ is also orthogonally complemented in $H$
such that $P_{\overline{\mathcal{R}(QP)}}=Q-P_{\mathcal{R}(Q)\cap \mathcal{N}(P)}$.
\end{lemma}

Next, we provide a useful lemma as follows.
\begin{lemma}\quad For every $P,Q\in\mathcal{P}(H)$, we have
\begin{equation}\label{relation T12 and T12*}
|P(I-Q)|+|(I-P)Q|=|(I-Q)P|+|Q(I-P)|.
\end{equation}
\end{lemma}
\begin{proof}[\bf Proof]For simplicity, we put
\begin{equation}\label{definition of T1 and T2} T_1=P(I-Q)  \quad  \mbox{and} \quad  T_2=(I-P)Q.
\end{equation}
It is clear that
$$T_1^*T_1=(I-Q)T_1^*T_1(I-Q)\quad\mbox{and}\quad
T_2^*T_2=QT_2^*T_2Q,$$
so
$$|T_1|=(I-Q)|T_1|(I-Q)\quad\mbox{and}\quad
|T_2|=Q|T_2|Q,$$
hence $|T_1|\cdot |T_2|=0$. Similarly, $|T_1^*|\cdot |T_2^*|=0$. As a result,
\begin{align*}
\big(|T_1|+|T_2|\big)^2=&T_1^*T_1+T_2^*T_2=P+Q-PQ-QP=T_1T_1^*+T_2T_2^*\\
=&\big(|T_1^*|+|T_2^*|\big)^2,
\end{align*}
which gives  \eqref{relation T12 and T12*} by taking the square roots of positive operators.
\end{proof}

Now, we are in the position to provide the main result of this section.
\begin{theorem}\quad  \label{technical thm} Let $P,Q\in\mathcal{P}(H)$ be such that $(P,Q)$ is semi-harmonious. Then there exists a unitary $U\in \mathcal{L}(H)$ such that
\begin{align}\label{equation of u equivalence-1}
&U(Q-P_{\mathcal{R}}) U^*=I-P-P_{\mathcal{N}},\\
\label{equation of u equivalence-2}&U(P-P_{\mathcal{R}}) U^*=I-Q-P_{\mathcal{N}},
\end{align}
where $\mathcal{R}$ and $\mathcal{N}$ are defined by \eqref{defn of R and N}.
\end{theorem}
\begin{proof}[\bf Proof] Let $T_1$ and $T_2$ be defined by \eqref{definition of T1 and T2}.
Since  $\overline{\mathcal{R}(P+Q)}$ is orthogonally complemented in $H$, by Lemma~\ref{charac of orth-PQ} both $\overline{\mathcal{R}(T_1^*)}$ and $\overline{\mathcal{R}(T_2)}$ are orthogonally complemented in $H$. Similarly, the orthogonal complementarity of $\overline{\mathcal{R}(2I-P-Q)}$ leads to that of $\overline{\mathcal{R}(T_1)}$ and $\overline{\mathcal{R}(T_2^*)}$. So for $i=1,2$, the notations $P_{\overline{\mathcal{R}(T_i)}}$ and $P_{\overline{\mathcal{R}(T^*_i)}}$ are meaningful. The point is, these projections can be used to obtain the canonical forms of $T_1$ and $T_2$ \cite{Corach-Maestripieri}. In fact, in view of \eqref{definition of T1 and T2} we have
\begin{equation*}
\overline{\mathcal{R}(T_1)}\subseteq \mathcal{R}(P)\quad \mbox{and}\quad\overline{\mathcal{R}(T^*_1)}\subseteq \mathcal{R}(I-Q),
\end{equation*}
hence
\begin{equation*}
P_{\overline{\mathcal{R}(T_1)}}P=P_{\overline{\mathcal{R}(T_1)}}
\quad \mbox{and} \quad
(I-Q)P_{\overline{\mathcal{R}(T^*_1)}}=P_{\overline{\mathcal{R}(T^*_1)}},
\end{equation*}
which lead to
\begin{equation}\label{alternative expression of T1}
T_1=P_{\overline{\mathcal{R}(T_1)}}T_1 P_{\overline{\mathcal{R}(T^*_1)}}=P_{\overline{\mathcal{R}(T_1)}}P_{\overline{\mathcal{R}(T^*_1)}}.
\end{equation}
Replacing $P$ and $Q$ with $I-P$ and $I-Q$ respectively, we obtain
\begin{equation}\label{alternative expression of T2}
T_2
=P_{\overline{\mathcal{R}(T_2)}}P_{\overline{\mathcal{R}(T^*_2)}}.
\end{equation}
Taking $*$-operation, from \eqref{alternative expression of T1} and \eqref{alternative expression of T2} we arrive at
\begin{equation}\label{alternative expression of T1* and  T2*}
T^*_1
=P_{\overline{\mathcal{R}(T^*_1)}}P_{\overline{\mathcal{R}(T_1)}}\quad\mbox{and}\quad  T^*_2
=P_{\overline{\mathcal{R}(T^*_2)}}P_{\overline{\mathcal{R}(T_2)}}.
\end{equation}

To make use of \eqref{relation T12 and T12*}, we need the polar decompositions of $T_1$ and $T_2$.
For $i=1,2$,  as $\overline{\mathcal{R}(T_i)}$ and $\overline{\mathcal{R}(T^*_i)}$ are both orthogonally complemented in $H$,
by \cite[Lemma~3.6 and Theorem~3.8]{Liu-Luo-Xu}  there exists a unique partial isometry  $V_i\in \mathcal{L}(H)$ such that
\begin{equation}\label{polar decomposition of Ti}
T_i=V_i |T_i|, \quad T^*_i=V^*_i |T^*_i|, \quad V^*_iV_i=P_{\overline{\mathcal{R}(T^*_i)}},  \quad V_i V^*_i=P_{\overline{\mathcal{R}(T_i)}}.
\end{equation}
Combining the last two equations in \eqref{polar decomposition of Ti} with \eqref{alternative expression of T1}--\eqref{alternative expression of T1* and  T2*}, we obtain
$$T_i=V_i V^*_iV^*_iV_i \quad \mbox{and}\quad T^*_i=V^*_iV_iV_iV^*_i.$$
This two equations together with \eqref{polar decomposition of Ti}, Lemmas~\ref{lem:Range Closure of T alpha and T} and \ref{lem:Range closure of TT and T} yield
\begin{equation*}
(V^*_i)^2 V_i=V_i^*(V_i V^*_iV^*_iV_i)=V^*_iT_i=V_i^*V_i |T_i|=|T_i|,
\end{equation*}
which gives
$$(V^*_i)^2 V_i=|T_i|=V_i^*V_i^2$$
by taking $*$-operation. Similarly, we have
\begin{equation*}\label{exchange star-1}V_i^2 V_i^*=|T_i^*|=V_i(V_i^*)^2.\end{equation*}
Consequently, equation~\eqref{relation T12 and T12*} turns out to be
\begin{equation}\label{equality of v12 and v12*}
V^*_1V_1^2+V^*_2V_2^2=V_1^2 V^*_1+V_2^2 V^*_2.
\end{equation}

Now, we are ready to construct the desired unitary operator.
Let
\begin{equation}\label{def of U12}U_1=V_1+P_{\mathcal{R}}, \quad U_2=V_2+P_{\mathcal{N}},\quad U=U_1-U_2,
\end{equation}
where $\mathcal{R}$ and $\mathcal{N}$ are defined by \eqref{defn of R and N}.
Then by Lemma~\ref{charac of orth-PQ}, \eqref{definition of T1 and T2} and \eqref{polar decomposition of Ti}, we have
\begin{align}\label{alter expression of p}
&P=V_1V^*_1+P_{\mathcal{R}},\quad   I-P=V_2V^*_2+P_{\mathcal{N}},\\
&\label{alter expression of Q}Q=V^*_2V_2+P_{\mathcal{R}},\quad  I-Q=V^*_1 V_1+P_{\mathcal{N}}.
\end{align}
It follows from \eqref{alter expression of p} that
$$V_1V_1^*P_{\mathcal{R}}=V_2V_2^*P_{\mathcal{N}}=V_1V_1^*V_2V_2^*=V_1V_1^*P_{\mathcal{N}}=P_{\mathcal{R}}V_2V_2^*=P_{\mathcal{R}}P_{\mathcal{N}}=0,$$
or equivalently,
\begin{equation}\label{six terms zero-01}V_1^*P_{\mathcal{R}}=V_2^*P_{\mathcal{N}}=V_1^*V_2=V_1^*P_{\mathcal{N}}=P_{\mathcal{R}}V_2=P_{\mathcal{R}}P_{\mathcal{N}}=0.\end{equation}
Similarly, it can be inferred from \eqref{alter expression of Q} that
\begin{equation}\label{six terms zero-02}V_2P_{\mathcal{R}}=V_1P_{\mathcal{N}}=V_2V_1^*=V_2P_{\mathcal{N}}=P_{\mathcal{R}}V_1^*=0.\end{equation}
It follows from \eqref{def of U12}--\eqref{six terms zero-02} that
\begin{align*}
&U_1U^*_2=U^*_1U_2=0, \quad U_1U_1^*=P, \quad U_2U_2^*=I-P,\\
&U_1^* U_1+U_2^* U_2=V_1^*V_1+P_{\mathcal{R}}+V_2^*V_2+P_{\mathcal{N}}=Q+I-Q=I.
\end{align*}
Therefore
$UU^*=U^*U=I$, so the operator $U$ defined by \eqref{def of U12} is a unitary.

Finally, we check the validity of \eqref{equation of u equivalence-1} and \eqref{equation of u equivalence-2}. According to \eqref{alter expression of p} and \eqref{alter expression of Q}, we have
\begin{align*}
&Q-P_{\mathcal{R}}=V_2^*V_2, \quad I-P-P_{\mathcal{N}}=V_2V_2^*,\\
&P-P_{\mathcal{R}}=V_1V_1^*, \quad I-Q-P_{\mathcal{N}}=V_1^*V_1,
\end{align*}
and thus
\begin{equation*}V_1V_1^*+V_2V_2^*=I-P_{\mathcal{R}}-P_{\mathcal{N}}=V_1^*V_1+V_2^*V_2.
\end{equation*}
The above equations together with \eqref{def of U12}, \eqref{six terms zero-01}, \eqref{six terms zero-02} and  \eqref{equality of v12 and v12*}  yield
\begin{align*}
U(Q-P_{\mathcal{R}})=&(U_1-U_2)V^*_2V_2=(V_1+P_{\mathcal{R}}-V_2-P_{\mathcal{N}})V^*_2V_2=-V_2\\
=&V_2V_2^*(V_1+P_{\mathcal{R}}-V_2-P_{\mathcal{N}})=(I-P-P_{\mathcal{N}})U,\\
U(P-P_{\mathcal{R}})=&(V_1+P_{\mathcal{R}}-V_2-P_{\mathcal{N}})V_1V_1^*=V_1^2V_1^*-V_2V_1V_1^*\\
=&V_1^2V_1^*-V_2(I-P_{\mathcal{R}}-P_{\mathcal{N}}-V_2V_2^*)=V_1^2V_1^*+V_2^2V_2^*-V_2\\
=&V^*_1V_1^2+V^*_2V_2^2-V_2=V^*_1V_1^2+(I-P_{\mathcal{R}}-P_{\mathcal{N}}-V_1^*V_1)V_2-V_2\\
=&V_1^*V_1^2-V_1^*V_1V_2=V_1^*V_1(V_1+P_{\mathcal{R}}-V_2-P_{\mathcal{N}})\\
=&(I-Q-P_{\mathcal{N}})U.
\end{align*}
Therefore, \eqref{equation of u equivalence-1} and \eqref{equation of u equivalence-2} are satisfied.
\end{proof}

\begin{remark}\quad\label{rem:harmonious case}{\rm  Let $P,Q\in\mathcal{P}(H)$ be such that $(P,Q)$ is harmonious. In this case, the Halmos' two projections theorem \cite[Theorem~3.3]{Xu-Yan} indicates that up to unitary equivalence, $P$ and $Q$ have the block matrix forms
\begin{equation*}\label{equ:block forms of P and Q}P=\left(
                  \begin{array}{cccccc}
                    I_{H_1} &  &  &  &  &  \\
                     & I_{H_2} &  &  &  &  \\
                     &  & 0 & &  &  \\
                     &  &  &  0 &  &  \\
                     &  &  &  & I_{H_5}&  \\
                     &  &  &  &  & 0 \\
                  \end{array}
                \right), \quad Q=\left(
                  \begin{array}{ccccc}
                    I_{H_1} &  &  &  &   \\
                     & 0 &  &  &    \\
                     &  & I_{H_3} & &   \\
                     &  &  &  0 &    \\
                     &  &  &  & Q_0  \\
                    \end{array}
                \right),\\
\end{equation*}
where
\begin{align}&H_1=\mathcal{R},\ H_2=\mathcal{R}(P)\cap\mathcal{N}(Q), \ H_3=\mathcal{N}(P)\cap\mathcal{R}(Q), \ H_4=\mathcal{N},\nonumber\\
 &H_5=\mathcal{R}(P)\ominus \left(H_1\oplus H_2\right),\ H_6=\mathcal{N}(P)\ominus \left(H_3\oplus H_4\right),\nonumber\\
\label{form of T}&Q_0=\left(
                    \begin{array}{cc}
                      A & A^\frac12\big(I_{H_5}-A\big)^\frac12 U_0^* \\
                     U_0A^\frac12\big(I_{H_5}-A\big)^\frac12  & U_0\big(I_{H_5}-A\big)U_0^*\\
                    \end{array}
                  \right)\in\mathcal{L}(H_5\oplus H_6),
\end{align}
in which $U_0\in\mathcal{L}(H_5,H_6)$ is a unitary, $A\in\mathcal{L}(H_5)$ is a positive contraction such that both $A$ and $I_{H_5}-A$ are injective and
$\overline{\mathcal{R}(A-A^2)}=H_5$, which implies that
$$\overline{\mathcal{R}(A)}=\overline{\mathcal{R}(I_{H_5}-A)}=H_5.$$
 With the notations as above and that in the proof of Theorem~\ref{technical thm}, we have
\begin{align*}&P_{\mathcal{R}}=I_{H_1}\oplus 0\oplus 0\oplus 0\oplus 0\oplus 0, \\
&P_{\mathcal{N}}=0\oplus 0\oplus 0 \oplus I_{H_4}\oplus 0\oplus 0,\\
&T_1=0\oplus I_{H_2}\oplus 0\oplus 0 \oplus S_1,\ V_1=0\oplus I_{H_2}\oplus 0\oplus 0 \oplus S_2,\\
&T_2=0\oplus 0 \oplus I_{H_3} \oplus 0 \oplus S_3,\ V_2=0\oplus 0 \oplus I_{H_3} \oplus 0 \oplus S_4,\\
\end{align*}
where
\begin{align*}&S_1=\left(
                    \begin{array}{cc}
                      I_{H_5}-A & -A^\frac12\big(I_{H_5}-A\big)^\frac12 U_0^* \\
                     0 & 0\\
                    \end{array}
                  \right),\\
&S_2=\left(
                    \begin{array}{cc}
                      \big(I_{H_5}-A\big)^\frac12  &-A^\frac12 U_0^* \\
                     0 & 0\\
                    \end{array}
                  \right),\\
&S_3=\left(
                    \begin{array}{cc}
                      0 & 0 \\
                     U_0A^\frac12\big(I_{H_5}-A\big)^\frac12  & U_0\big(I_{H_5}-A\big)U_0^*\\
                    \end{array}
                  \right),\\
&S_4=\left(
                    \begin{array}{cc}
                      0 & 0 \\
                     U_0A^\frac12 & U_0\big(I_{H_5}-A\big)^\frac12 U_0^*\\
                    \end{array}
                  \right).
\end{align*}
In virtue of  \eqref{def of U12}, we have
\begin{equation*}U=V_1+P_{\mathcal{R}}-V_2-P_{\mathcal{N}}=\mbox{diag}(X,Y),
\end{equation*}
where $X=\mbox{diag}(I_{H_1}, I_{H_2}, -I_{H_3}, -I_{H_4})$ and
\begin{equation*}Y=S_2-S_4=\left(
    \begin{array}{cc}
        \big(I_{H_5}-A\big)^\frac12  &-A^\frac12 U_0^*  \\
      -U_0A^\frac12 & -U_0\big(I_{H_5}-A\big)^\frac12 U_0^* \\
     \end{array}
  \right).
\end{equation*}
This gives the block matrix form of the unitary operator $U$ satisfying \eqref{equation of u equivalence-1} and \eqref{equation of u equivalence-2}.
}\end{remark}

\begin{remark}\quad{\rm
Since every closed linear subspace of a Hilbert space is orthogonally complemented,
Theorem~\ref{technical thm}  is therefore always applicable to every pair of projections on a Hilbert space.
}\end{remark}

\begin{theorem}\quad  \label{thm:unitary equivalence stops}There exist projections $P$ and $Q$ on certain Hilbert $C^*$-module such that $\mathcal{R}=\mathcal{N}=\{0\}$, whereas  \eqref{equation of u equivalence-1} and \eqref{equation of u equivalence-2} have no common unitary operator solution.
\end{theorem}
\begin{proof}[\bf Proof] Following \cite[Section~3]{Manuilov-Moslehian-Xu}, we put
  $\mathfrak{B}=C\big([0,1];M_2(\mathbb{C})\big)$ and set
\begin{equation}\label{A}
\mathfrak{A}=\{f\in \mathfrak{B}:f(0)\ \mbox{and}\ f(1)\mbox{\ are both\ diagonal}\}.
\end{equation}
As is shown in the proof of Theorem~\ref{thm:semi-harmonious not implication}, $\mathfrak{A}$ itself is a Hilbert $\mathfrak{A}$-module, and we can identify $\mathcal{L}(\mathfrak{A})$ with $\mathfrak{A}$. Let $H=\mathfrak{A}$ and $Q\in\mathcal{P}(H)$ be determined by \eqref{equ:defn of P t and Q t}, and let $P\in\mathcal{P}(H)$ be changed to
\begin{equation*}\label{new defn of P (t)}P(t)=\left(\begin{matrix}c_t^2&-s_tc_t\\-s_tc_t&s_t^2\end{matrix}\right), \quad \mbox{for $t\in [0,1]$},
\end{equation*}
where $c_t$ and $s_t$ are defined by \eqref{defn of s t and c t}.

Let $\mathcal{R}$ and $\mathcal{N}$ be defined by \eqref{defn of R and N}, and suppose that $x\in \mathcal{N}$ is determined by $x(t)=\big(x_{ij}(t)\big)_{1\le i,j\le 2}$ for $t\in [0,1]$. Utilizing
\begin{equation}\label{equ:expression of P plus Q t}P(t)+Q(t)=\left(
             \begin{array}{cc}
               2c_t^2 & 0 \\
               0 & 2s_t^2 \\
             \end{array}
           \right)\quad\mbox{and}\quad \big[P(t)+Q(t)\big]x(t)=0
           \end{equation} for $t\in [0,1]$, we obtain $x_{ij}(t)=0$ for $i,j\in \{1,2\}$ and  every $t\in (0,1)$,
which imply that  $x_{ij}=0$ for $1\le i,j\le 2$,  since all functions considered are continuous on $[0,1]$.  This shows that $\mathcal{N}=\{0\}$.
In view of $\mathcal{R}=\mathcal{N}(I-P)\cap \mathcal{N}(I-Q)$, the proof of  $\mathcal{R}=\{0\}$ is similar.

Suppose that $U$ determined by $U(t)=\big(U_{ij}(t)\big)_{1\le i,j\le 2}$ is a unitary in $H$, which satisfies both \eqref{equation of u equivalence-1} and \eqref{equation of u equivalence-2}. Due to $P_{\mathcal{R}}=0$ and $P_{\mathcal{N}}=0$, \eqref{equation of u equivalence-1} and \eqref{equation of u equivalence-2} are simplified as
$$UP=(I-Q)U\quad\mbox{and}\quad UQ=(I-P)U,$$
or equivalently,
\begin{equation}\label{equ:two derived equations}U(Q-P)=(Q-P)U\quad\mbox{and}\quad U(P+Q)=(2I-P-Q)U.
\end{equation}
Substituting
$$(Q-P)(t)=\left(
                 \begin{array}{cc}
                   0 & 2s_tc_t \\
                   2s_tc_t & 0 \\
                 \end{array}
               \right), \quad t\in [0,1],
               $$
into the first equation in \eqref{equ:two derived equations} yields
$$U_{12}(t)=U_{21}(t), \quad U_{11}(t)=U_{22}(t), \quad t\in [0,1].$$
Combining the above equations with the expression of $P(t)+Q(t)$ given in \eqref{equ:expression of P plus Q t} and
the second equation in \eqref{equ:two derived equations}, we arrive at
$$U_{11}(t)\left(c_t^2-s_t^2\right)=0,\quad \forall t\in [0,1],$$ hence
$U_{11}(t)\equiv 0$ for $t\in [0,1]$ by the continuity of $U_{11}$. Consequently,
 $$U(t)=\left(
                                                                   \begin{array}{cc}
                                                                     0 & U_{12}(t) \\
                                                                     U_{12}(t) & 0 \\
                                                                   \end{array}
                                                                 \right).$$ This together with \eqref{A} yields
$U(0)=U(1)=\left(
        \begin{array}{cc}
          0 & 0 \\
          0 & 0 \\
        \end{array}
      \right)$. It is a contradiction, since $U$ is a unitary in $H$ which ensures that all the $2\times 2$ matrices $U(t) (t\in [0,1])$  are unitary.
 \end{proof}

 \section{$C^*$-isomorphisms associated with two projections}\label{sec:norm equations}
Unless otherwise specified, throughout this section $(P,Q,H)$ is a matched triple, $i(P,Q,H)$ and $o(P,Q,H)$ are its inner algebra and  outer algebra (see Definitions~\ref{defn of matched triple}).
It is clear that
\begin{equation}\label{dense subset}i(P,Q,H)=\overline{\mbox{span}\left\{X^{(P,Q,k)}: X\in\{A,B,C,D\}, k\in\mathbb{Z}_+\right\}},\end{equation}
where $\mathcal{R}$ (also $\mathcal{N}$) is defined by \eqref{defn of R and N}, $P_{\mathcal{R}}=P\wedge Q$ and
\begin{align}\label{defn of A 1}&A^{(P,Q,k)}=\left[(P-P_{\mathcal{R}})(Q-P_{\mathcal{R}})\right]^k,\\
 \label{defn of B 1}&B^{(P,Q,k)}=A^{(P,Q,k)}(P-P_{\mathcal{R}}),\\
\label{defn of C 1} &C^{(P,Q,k)}=\left(A^{(P,Q,k)}\right)^*=A^{(Q,P,k)},\\
 \label{defn of D 1} &D^{(P,Q,k)}=C^{(P,Q,k)}(Q-P_{\mathcal{R}})=B^{(Q,P,k)}
 \end{align}
with the convention that $A^{(P,Q,0)}=I$. For each $k\ge 1$, by utilizing
$P_{\mathcal{R}}\le P$ and $P_{\mathcal{R}}\le Q$ we obtain
\begin{align}\label{A power k}A^{(P,Q,k)}=&(PQ-P_{\mathcal{R}})^k=(PQ)^k-P_{\mathcal{R}},
\end{align}
which gives
\begin{align}A^{(P,Q,k)}\left(A^{(P,Q,k)}\right)^*=&\left[(PQ)^k-P_{\mathcal{R}}\right]\left[(QP)^k-P_{\mathcal{R}}\right]\nonumber\\
\label{computations of power-01}=&(PQ)^k(QP)^k-P_{\mathcal{R}}=(PQP)^{2k-1}-P_{\mathcal{R}}\\
=&\left(PQP-P_{\mathcal{R}}\right)^{2k-1}.\nonumber
\end{align}
It follows that
\begin{equation}\label{computations of power-02}A^{(P,Q,k)}\left(A^{(P,Q,k)}\right)^*=\left[A^{(P,Q,1)}\left(A^{(P,Q,1)}\right)^*\right]^{2k-1},\quad \forall k\ge 1.\end{equation}
Similarly,
\begin{align}\label{computations of power-03}B^{(P,Q,k)}=&\left[(PQ)^k-P_{\mathcal{R}}\right](P-P_{\mathcal{R}})=(PQP)^k-P_{\mathcal{R}}=\left(PQP-P_{\mathcal{R}}\right)^k\nonumber\\
=&\left[A^{(P,Q,1)}\left(A^{(P,Q,1)}\right)^*\right]^k, \quad \forall k\ge 1.
\end{align}

Now, let $(\pi, X)$ be a faithful representation of $o(P,Q,H)$.  Replacing $X$ with $\pi(I)X$ if necessary, in what follows we always assume that $\pi$ is unital.
The infimum of $\pi(P)$ and $\pi(Q)$, and its subtraction by
$\pi(P_{\mathcal{R}})$ are denoted simply by $P_{\mathcal{R}\pi}$ and $\widetilde{P_{\mathcal{R}}}$ respectively, that is,
\begin{equation}\label{p pi R}P_{\mathcal{R}\pi}=P_{\mathcal{R}\left(\pi(P)\right)\cap \mathcal{R}\left(\pi(Q)\right)},\quad
\widetilde{P_{\mathcal{R}}}=P_{\mathcal{R}\pi}-\pi(P_{\mathcal{R}}).
\end{equation}
It is clear that $\pi(P_{\mathcal{R}})\le P_{\mathcal{R}\pi}$, so $\widetilde{P_{\mathcal{R}}}$ defined as above is a projection. By \eqref{defn of A 1} we have $A^{(\pi(P),\pi(Q),0)}=I$ and when $k\ge 1$,
\begin{align*}A^{(\pi(P),\pi(Q),k)}=&\left[(\pi(P)-P_{\mathcal{R}\pi})(\pi(Q)-P_{\mathcal{R}\pi})\right]^k=\left[\pi(P)\pi(Q)\right]^k-P_{\mathcal{R}\pi}.
\end{align*}
This together with \eqref{A power k} and \eqref{p pi R} gives
$$\pi\left(A^{(P,Q,k)}\right)=A^{(\pi(P),\pi(Q),k)}+\widetilde{P_{\mathcal{R}}}.$$
The derivation above shows that
\begin{equation}\label{ley relation for norm-1}\pi\left(X^{(P,Q,k)}\right)=X^{(\pi(P),\pi(Q),k)}+\widetilde{P_{\mathcal{R}}}\quad\mbox{whenever $X^{(P,Q,k)}\ne I$}.
\end{equation}
Note that $\widetilde{P_{\mathcal{R}}}$ is a projection,   and
\begin{equation}\label{ley relation for norm-111}X^{\left(\pi(P),\pi(Q),k\right)}\cdot \widetilde{P_{\mathcal{R}}}=\widetilde{P_{\mathcal{R}}}\cdot X^{\left(\pi(P),\pi(Q),k\right)}=0\quad\mbox{whenever $X^{(P,Q,k)}\ne I$},\end{equation}
so by \eqref{ley relation for norm-1} together with the observation $\|I\|=1\ge \|\widetilde{P_{\mathcal{R}}}\|$, we arrive at
\begin{align}\left\|X^{(P,Q,k)}\right\|=&\left\|\pi\left(X^{(P,Q,k)}\right)\right\|\nonumber\\
=&\max\left\{\left\|X^{(\pi(P),\pi(Q),k)}\right\|, \, \|\widetilde{P_{\mathcal{R}}}\|\right\}, \quad \forall k\ge 0.
\end{align}

We are now ready to derive a couple of norm equations. The first one reads as follows.

\begin{lemma}\quad\label{technical lem1-norm} Let $(\pi, X)$ be a faithful representation of $o(P,Q,H)$. Then for every $X\in\{A,B,C,D\}$ and $k\in\mathbb{Z}_+$, we have
\begin{equation}\label{equ:fu-xu-00}\left\|X^{(P,Q,k)}\right\|=\left\|X^{(\pi(P),\pi(Q),k)}\right\|.\end{equation}
 \end{lemma}
\begin{proof}[\bf Proof]First, we prove that
\begin{equation}\label{equ:fu-xu-01}\left\|A^{(P,Q,k)}\right\|=\left\|A^{(\pi(P),\pi(Q),k)}\right\|, \quad \forall k\in\mathbb{Z}_+.
\end{equation}
The case of $k=0$ is trivial, so we start with $k=1$. From \eqref{ley relation for norm-1}, we have
\begin{align}\left\|A^{(P,Q,1)}\right\|^2=&\left\|\pi\left(A^{(P,Q,1)}\right)\right\|^2=\left\|\pi\left(A^{(P,Q,1)}\right)\left[\pi\left(A^{(P,Q,1)}\right)\right]^*\right\|\nonumber\\
=&\left\|\left(A^{(\pi(P),\pi(Q),1)}+\widetilde{P_{\mathcal{R}}}\right)\left[\left(A^{(\pi(P),\pi(Q),1)}\right)^*+\widetilde{P_{\mathcal{R}}}\right]\right\|\nonumber\\
=&\left\|A^{(\pi(P),\pi(Q),1)}\left(A^{(\pi(P),\pi(Q),1)}\right)^*+\widetilde{P_{\mathcal{R}}}\right\|\nonumber\\
=&\max\left\{\left\|A^{(\pi(P),\pi(Q),1)}\left(A^{(\pi(P),\pi(Q),1)}\right)^*\right\|,\, \|\widetilde{P_{\mathcal{R}}}\|\right\}\nonumber\\
\label{computations of power-022}=&\max\left\{\left\|A^{(\pi(P),\pi(Q),1)}\right\|^2, \, \|\widetilde{P_{\mathcal{R}}}\|\right\},
\end{align}
which implies that
\begin{equation}\label{equ:fu-xu-01-special}\left\|A^{(P,Q,1)}\right\|=\left\|A^{(\pi(P),\pi(Q),1)}\right\|\end{equation}
whenever $\left\|A^{(\pi(P),\pi(Q),1)}\right\|=1$. Suppose that $\left\|A^{(\pi(P),\pi(Q),1)}\right\|<1$. Then according to \eqref{computations of power-01} and \eqref{computations of power-02}, we have
\begin{align*}\left\|\big(\pi(P)\pi(Q)\pi(P)\big)^{2k-1}-P_{\mathcal{R}\pi}\right\|=&\left\|\left[A^{(\pi(P),\pi(Q),1)}\left(A^{(\pi(P),\pi(Q),1)}\right)^*\right]^{2k-1}\right\|\\
=&\left\|A^{(\pi(P),\pi(Q),1)}\left(A^{(\pi(P),\pi(Q),1)}\right)^*\right\|^{2k-1}\\
=&\left\|A^{(\pi(P),\pi(Q),1)}\right\|^{4k-2}\to 0\quad \mbox{as $k\to \infty$.}
\end{align*}
It follows that $P_{\mathcal{R}\pi}$ is contained in the $C^*$-algebra $\pi(C^*_{(P,Q)})$, so
$\pi^{-1}(P_{\mathcal{R}\pi})\le \pi^{-1}\big(\pi(P)\big)=P$ and $\pi^{-1}(P_{\mathcal{R}\pi})\le Q$ as well. Therefore, $\pi^{-1}(P_{\mathcal{R}\pi})\le P_{\mathcal{R}}$, and thus $P_{\mathcal{R}\pi}\le \pi(P_\mathcal{R})\le P_{\mathcal{R}\pi}$. Plugging
$\widetilde{P_{\mathcal{R}}}=0$ into  \eqref{computations of power-022} gives \eqref{equ:fu-xu-01-special} immediately.

Now, we consider the case that $k\ge 2$.  Due to \eqref{computations of power-02} and  \eqref{equ:fu-xu-01-special}, we have
\begin{align*}\left\|A^{(P,Q,k)}\right\|=\left\|A^{(P,Q,1)}\right\|^{2k-1}=\left\|A^{(\pi(P),\pi(Q),1)}\right\|^{2k-1}
=\left\|A^{(\pi(P),\pi(Q),k)}\right\|.
\end{align*}
This completes the proof of
\eqref{equ:fu-xu-01}.

Next, we prove that
\begin{equation}\label{equ:fu-xu-02}\left\|B^{(P,Q,k)}\right\|=\left\|B^{(\pi(P),\pi(Q),k)}\right\|, \quad \forall k\in\mathbb{Z}_+.
\end{equation}
By \eqref{defn of B 1}  both $B^{(P,Q,0)}$ and $B^{(\pi(P),\pi(Q),0)}$ are projections, and
\begin{equation*}\left\|B^{(P,Q,0)}\right\|=0\Longleftrightarrow P\le Q\Longleftrightarrow \pi(P)\le \pi(Q)\Longleftrightarrow \left\|B^{(\pi(P),\pi(Q),0)}\right\|=0.
\end{equation*}
This shows the validity of \eqref{equ:fu-xu-02} for $k=0$. Suppose that $k\ge 1$. Then by \eqref{computations of power-03}, we can get
$$\left\|B^{(P,Q,k)}\right\|=\left\|A^{(P,Q,1)}\right\|^{2k}=\left\|A^{(\pi(P),\pi(Q),1)}\right\|^{2k}=\left\|B^{(\pi(P),\pi(Q),k)}\right\|.$$

Exchanging $P$ with $Q$, we conclude that for every $k\in\mathbb{Z}_+$,
\begin{align*}&\left\|C^{(P,Q,k)}\right\|=\left\|A^{(Q,P,k)}\right\|=\left\|A^{(\pi(Q),\pi(P),k)}\right\|=\left\|C^{(\pi(P),\pi(Q),k)}\right\|,\\
&\left\|D^{(P,Q,k)}\right\|=\left\|B^{(Q,P,k)}\right\|=\left\|B^{(\pi(Q),\pi(P),k)}\right\|=\left\|D^{(\pi(P),\pi(Q),k)}\right\|.
\end{align*}
This completes the proof of \eqref{equ:fu-xu-00}.
\end{proof}

\begin{corollary}\quad Let $(\pi, X)$ be a faithful representation of $o(P,Q,H)$. Then for every $n\in\mathbb{N}$, $X_i\in\{A,B,C,D\}$ and $k_i\in\mathbb{Z}_+$  $(1\le i\le n)$,
we have
\begin{equation}\label{norm equation finite product}\left\|\prod_{i=1}^n X_i^{(P,Q,k_i)}\right\|=\left\|\prod_{i=1}^n X_i^{(\pi(P),\pi(Q),k_i)}\right\|.\end{equation}
 \end{corollary}
\begin{proof}[\bf Proof]From the definition of $X^{(P,Q,k)}$ given by \eqref{defn of A 1}--\eqref{defn of D 1}, it is clear that
$$\prod_{i=1}^n X_i^{(P,Q,k_i)}=Z^{(P,Q,k)}\quad\mbox{and}\quad \prod_{i=1}^n X_i^{(\pi(P),\pi(Q),k_i)}=Z^{(\pi(P),\pi(Q),k)}$$
for some $Z\in\{A,B,C,D\}$ and $k\in\mathbb{Z}_+$. Due to \eqref{equ:fu-xu-00}, the desired norm equation follows.
\end{proof}

We provide a technical lemma as follows.

 \begin{lemma}\quad\label{technical lem2-norm} Let $(P,Q)$ be a harmonious pair of projections on $H$. Suppose that $n\in\mathbb{N}$, $X_i\in\{A,B,C,D\}$ and $k_i\in\mathbb{Z}_+$ $(1\le i\le n)$ are given such that \begin{equation}\label{equ: norm equals one}\left\|\prod_{i=1}^n X_i^{(P,Q,k_i)}\right\|=1.\end{equation} Then for every $\lambda_i\in \mathbb{C}$  $(1\le i\le n)$, we have
\begin{equation}\label{norm lower bound-01}\left|\sum_{i=1}^n \lambda_i\right|\le \left\|\sum_{i=1}^n \lambda_i X_i^{(P,Q,k_i)}\right\|.
\end{equation}
\end{lemma}
\begin{proof}[\bf Proof]Denote by $\lambda=\sum\limits_{i=1}^n \lambda_i$. The verification of
\begin{equation*}|\lambda|\le \left\|\sum_{i=1}^n \lambda_i X_i^{(P,Q,k_i)}\right\|
\end{equation*}
will be carried out by taking several cases into consideration.

\textbf{Case 1}\quad $X_i^{(P,Q,k_i)}\in \{I, P-P_{\mathcal{R}}\}$ for all $i\in\{1,2\dots,n\}$.
If $X_i^{(P,Q,k_i)}\equiv I$, then \eqref{norm lower bound-01} is obviously satisfied. Otherwise, we have
\begin{equation*}\label{becomes one operator-01}\prod_{i=1}^n X_i^{(P,Q,k_i)}=P-P_{\mathcal{R}}, \end{equation*}
so according to \eqref{equ: norm equals one} we obtain $\|P-P_{\mathcal{R}}\|=1$.
It follows that
\begin{align*}&\left\|\sum_{i=1}^n \lambda_i X_i^{(P,Q,k_i)}\right\|\ge \left\|(P-P_{\mathcal{R}})\sum_{i=1}^n \lambda_i X_i^{(P,Q,k_i)}\right\|=
\left\|\lambda (P-P_{\mathcal{R}})\right\|=|\lambda|.
\end{align*}

\textbf{Case 2}\quad $X_i^{(P,Q,k_i)}\in \{I, Q-P_{\mathcal{R}}\}$ for all $i\in\{1,2\dots,n\}$. The same verification gives \eqref{norm lower bound-01}.

\textbf{Case 3}\quad There exist $i_1,i_2\in\{1,2,\dots,n\}$ such that $X_{i_1}^{(P,Q,k_{i_1})}\notin \{I, P-P_{\mathcal{R}}\}$ and
$X_{i_2}^{(P,Q,k_{i_2})}\notin \{I, Q-P_{\mathcal{R}}\}$. In this case, firstly we  show  that
\begin{equation}\label{norm equals one single term} \left\|A^{(P,Q,1)}\right\|=1.\end{equation}

\textbf{Subcase 1}\quad $k_{i_1}\ne 0$ or $k_{i_2}\ne 0$. Without loss of generality, we may assume that $k_{i_1}\ne 0$. In this subcase,
$$1\ge \left\|A^{(P,Q,1)}\right\|=\left\|C^{(P,Q,1)}\right\|\ge \left\|X_{i_1}^{(P,Q,k_{i_1})}\right\|\ge \left\|\prod_{i=1}^n X_i^{(P,Q,k_i)}\right\|=1.$$
Thus, \eqref{norm equals one single term} is satisfied.

\textbf{Subcase 2}\quad $k_{i_1}=k_{i_2}=0$. In this subcase, we have $X_{i_1}^{(P,Q,k_{i_1})}=Q-P_{\mathcal{R}}$ and $X_{i_2}^{(P,Q,k_{i_2})}=P-P_{\mathcal{R}}$, which mean that
$$\prod_{i=1}^n X_i^{(P,Q,k_i)}=W_1 A^{(P,Q,1)}W_2\quad \mbox{or}\quad \prod_{i=1}^n X_i^{(P,Q,k_i)}=W_1 C^{(P,Q,1)}W_2$$ for some contractions $W_1,W_2\in \mathcal{L}(H)$. As is shown in \textbf{Subcase 1}, \eqref{norm equals one single term} is also satisfied.

 Since $(P,Q)$ is harmonious, the Halmos' two projections theorem is applicable. Following the notations as in Remark~\ref{rem:harmonious case}, we have\footnote{\,In some cases, it may happen that the closed subspaces  $H_5$ and $H_6$  constructed for the Halmos decomposition are trivial, that is, $H_5=H_6=\{0\}$. Due to \eqref{norm equals one single term}, both $H_5$ and $H_6$ are non-trivial, hence $A^{(P,Q,1)}$ has the form \eqref{H 5 appear}.}
\begin{align}&P-P_{\mathcal{R}}=0\oplus I_{H_2}\oplus 0\oplus 0\oplus I_{H_5}\oplus 0,\nonumber\\
&Q-P_{\mathcal{R}}=0\oplus 0\oplus I_{H_3}\oplus 0\oplus Q_0,\nonumber\\
\label{H 5 appear}&A^{(P,Q,1)}=0\oplus 0\oplus 0\oplus 0 \oplus S,
\end{align}
where $Q_0$ is defined by \eqref{form of T} and $S$ is given by
\begin{equation*}\label{form of S}S=\left(
                    \begin{array}{cc}
                      A & A^\frac12\big(I_{H_5}-A\big)^\frac12 U_0^* \\
                     0  & 0\\
                    \end{array}
                  \right).
\end{equation*}
Based on the above block matrices, we have
\begin{equation}\label{norm A equals one}\|A\|=\left\|\left(
                 \begin{array}{cc}
                   A & 0 \\
                   0 & 0 \\
                 \end{array}
               \right)\right\|=
\|SS^*\|=\|S\|^2=\left\|A^{(P,Q,1)}\right\|^2=1.\end{equation}
This together with the positivity and contraction of $A$ implies that $1\in \mbox{sp}(A)$, where $\mbox{sp}(A)$ denotes the spectrum of $A$.

It is easy to verify that for every $X\in\{A,B,C,D\}$ and $k\in\mathbb{Z}_+$, there exists $r\in \mathbb{N}$ depending on $X$ and $k$ such that
$$A^{(P,Q,1)}X^{(P,Q,k)}(P-P_{\mathcal{R}})=0\oplus 0\oplus 0\oplus 0\oplus \left(
                                                                              \begin{array}{cc}
                                                                                A^r &  \\
                                                                                 & 0 \\
                                                                              \end{array}
                                                                            \right).$$
Consequently,
\begin{align}\left\|\sum_{i=1}^n \lambda_i X_i^{(P,Q,k_i)}\right\|&\ge \left\|A^{(P,Q,1)}\left(\sum_{i=1}^n \lambda_i X_i^{(P,Q,k_i)}\right)(P-P_{\mathcal{R}})\right\|\nonumber\\
\label{norm lower bound-02}&=\left\|\left(
    \begin{array}{cc}
      \sum\limits_{i=1}^n \lambda_i A^{r_i} &  \\
       & 0\\
    \end{array}
  \right)\right\|=\left\|\sum\limits_{i=1}^n \lambda_i A^{r_i}\right\|
\end{align}
for some $r_i\in\mathbb{N}$ $(1\le i\le n)$. Let $f(t)=\sum\limits_{i=1}^n \lambda_i t^{r_i}$ for $t\ge 0$. Then
\begin{align*}\left\|\sum\limits_{i=1}^n \lambda_i A^{r_i}\right\|=\max\left\{|f(t)|:t\in \mbox{sp}(A)\right\}\ge |f(1)|=|\lambda|.
\end{align*}
Combining the above inequality with \eqref{norm lower bound-02} gives \eqref{norm lower bound-01}.
\end{proof}

Along the same line, another technical  lemma can be provided as follows.

 \begin{lemma}\quad\label{technical lem2-norm++} Let $(P,Q)$ be a harmonious pair of projections on $H$ such that $PQ\ne QP$. Suppose that $n\in\mathbb{N}$, $X_i\in\{A,B,C,D\}$ and $k_i\in\mathbb{Z}_+$ are given such that \eqref{equ: norm equals one} is satisfied and
 $X_i^{(P,Q,k_i)}\ne I$ for all $i\in\{1,2,\dots,n\}$. Then for every $\lambda_i\in \mathbb{C}$ $(0\le i\le n)$, we have
 \begin{equation}\label{norm lower bound-01++}|\lambda_0|\le \left\|\lambda_0 \left(I-P_{\mathcal{R}}\right)+\sum_{i=1}^n \lambda_i X_i^{(P,Q,k_i)}\right\|.
\end{equation}
\end{lemma}
\begin{proof}[\bf Proof]As in the proof of Lemma~\ref{technical lem2-norm}, the verification of \eqref{norm lower bound-01++} will be dealt with via several cases. Let $\lambda\in\mathbb{C}$ and $W\in\mathcal{L}(H)$ be defined by
\begin{equation}\label{defn of lambda and W}\lambda=\sum_{i=0}^n \lambda_i, \quad W=\lambda_0 \left(I-P_{\mathcal{R}}\right)+\sum_{i=1}^n \lambda_i X_i^{(P,Q,k_i)}.
\end{equation}

\textbf{Case 1}\quad $X_i^{(P,Q,k_i)}=P-P_{\mathcal{R}}$ for all $i\in\{1,2\dots,n\}$. In this case,
\begin{align*}W=\lambda_0 (I-P)+\lambda (P-P_{\mathcal{R}}).
\end{align*}
Since $PQ\ne QP$, we have $I-P\ne 0$, hence
$$\|W\|\ge \|\lambda_0 (I-P)\|=|\lambda_0|.$$

\textbf{Case 2}\quad $X_i^{(P,Q,k_i)}=Q-P_{\mathcal{R}}$ for all $i\in\{1,2\dots,n\}$. As is shown in the above case, we have
$\|W\|\ge |\lambda_0|$.

\textbf{Case 3}\quad There exist $i_1,i_2\in\{1,2,\dots,n\}$ such that $X_{i_1}^{(P,Q,k_{i_1})}\ne P-P_{\mathcal{R}}$ and
$X_{i_2}^{(P,Q,k_{i_2})}\ne Q-P_{\mathcal{R}}$.
 Following the notations as in Remark~\ref{rem:harmonious case}, we have $H=\bigoplus\limits_{i=1}^6H_i$ and up to unitary equivalence, every operator $Y\in\mathcal{L}(H)$ has the matrix form $Y=(Y_{ij})_{1\le i,j\le 6}$ with $Y_{ij}\in\mathcal{L}(H_j,H_i)$. Let the linear map $\phi:\mathcal{L}(H)\to \mathcal{L}(H_6)$ be defined by
 $\phi(Y)=Y_{66}$. According to the definition of  $X^{(P,Q,k)}$ given by \eqref{defn of A 1}--\eqref{defn of D 1}, we have
  $\phi\left(B^{(P,Q,0)}\right)=0$ and
$$\phi\left(A^{(P,Q,k)}\right)=\phi\left(B^{(P,Q,k)}\right)=\phi\left(C^{(P,Q,k)}\right)=0$$
for every $k\ge 1$. Furthermore, direct computations yield
$$\phi\left(D^{(P,Q,k)}\right)=U_0(I-A)A^kU_0^*,\quad \forall k\in\mathbb{Z}_+.$$
It follows from \eqref{defn of lambda and W} that\footnote{\,If $X_i\ne D$ for all $i\in\{1,2,\dots,n\}$, then $\phi(W)=\lambda_0 I$, so in this case $\|\phi(W)\|=|\lambda_0|$.}
\begin{align*}\phi(W)=U_0\Big[\lambda_0 I+\sum_j \lambda_{i_j} (I-A)A^{k_{i_j}}\Big]U_0^*,
\end{align*}
where $i_j$ is chosen in $\{1,2,\dots,n\}$ whenever $X_{i_j}^{(P,Q,k_{i_j})}=D^{(P,Q,k_{i_j})}$.
Let
$$f(t)=\lambda_0+\sum_j \lambda_{i_j} (1-t)t^{k_{i_j}}, \quad t\ge 0.$$
Since $0\le A\le I$ and $1\in \mbox{sp}(A)$ (see \eqref{norm A equals one}), we have
\begin{align*}\|W\|\ge \|\phi(W)\|=\max\left\{|f(t)|:t\in\mbox{sp}(A)\right\}\ge |f(1)|=|\lambda_0|.
\end{align*}
This completes the proof of \eqref{norm lower bound-01++}.
\end{proof}

Now, we provide the main result of this section as follows.

\begin{theorem}\quad  \label{thm:algebras isomorphic-01} For each faithful representation $(\pi,X)$  of $o(P,Q,H)$,  a faithful  representation $(\widetilde{\pi},X)$ of $i(P,Q,H)$  can be induced such that  $\widetilde{\pi}(I)=\pi(I)$, and
\begin{equation}\label{two natural map ranges}\widetilde{\pi}(P-P_{\mathcal{R}})=\pi(P)-\pi(P)\wedge \pi(Q),\quad \widetilde{\pi}(Q-P_{\mathcal{R}})=\pi(Q)-\pi(P)\wedge \pi(Q).
\end{equation}
\end{theorem}
\begin{proof}[\bf Proof]According to \eqref{dense subset}, we need only to prove that $\|T\|=\|\underline{T}\|$  for every  $n\in\mathbb{N}$, $X_i\in\{A,B,C,D\}$, $k_i\in\mathbb{Z}_+$ and  $\lambda_i\in\mathbb{C}$ ($0\le i\le n$) such that $X_i^{(P,Q,k_i)}\ne I$ for all $i\ge 1$,
where\footnote{\,When $\lambda_0=0$, the first term in $T$ and $\underline{T}$ will disappear.}
\begin{equation}\label{defn of T and tilde T}T=\lambda_0 I+\sum_{i=1}^n \lambda_i X_i^{(P,Q,k_i)}, \quad \underline{T}=\lambda_0I+\sum_{i=1}^n \lambda_i X_i^{(\pi(P),\pi(Q),k_i)}.
\end{equation}
If $\widetilde{P_{\mathcal{R}}}=0$, then we see from \eqref{ley relation for norm-1} that
$\pi(T)=\underline{T}$, which gives $\|T\|=\|\underline{T}\|$ as $\pi$ is faithful. In what follows, we assume that $\widetilde{P_{\mathcal{R}}}\ne 0$. In this case, we have $PQ\ne QP$. Otherwise, $P_{\mathcal{R}}=PQ$  and $P_{\mathcal{R}\pi}=\pi(P)\pi(Q)=\pi(P_{\mathcal{R}})$, which contradicts the assumption of $\widetilde{P_{\mathcal{R}}}\ne 0$.
Let  $\lambda=\sum\limits_{i=0}^n \lambda_i$.  By \eqref{ley relation for norm-1} we have
\begin{align}\pi(T)=\underline{T}+(\lambda-\lambda_0) \widetilde{P_{\mathcal{R}}}=L+\lambda \widetilde{P_{\mathcal{R}}},
\end{align}
where
\begin{equation*}L=\lambda_0\left(I-\widetilde{P_{\mathcal{R}}}\right)+\sum_{i=1}^n \lambda_i X_i^{(\pi(P),\pi(Q),k_i)}=Z+\lambda_0 \pi(P_{\mathcal{R}}),\end{equation*}
in which
$$Z=\lambda_0\left(I-P_{\mathcal{R}\pi}\right)+\sum_{i=1}^n \lambda_i X_i^{(\pi(P),\pi(Q),k_i)}.$$
Hence
\begin{align}\label{new estimation of L}\|L\|=&\max\{\|Z\|,\, |\lambda_0|\, \|\pi(P_{\mathcal{R}})\|\},\\
\label{new estimation of T}\|T\|=&\|\pi(T)\|=\max\left\{\|L\|, \, |\lambda|\, \|\widetilde{P_{\mathcal{R}}}\|\right\}=\max\left\{\|L\|, \, |\lambda|\right\},\\
\label{new estimation of T down}\|\underline{T}\|=&\|L+\lambda_0 \widetilde{P_{\mathcal{R}}}\|=\max\left\{\|L\|, \, |\lambda_0|\, \|\widetilde{P_{\mathcal{R}}}\|\right\}=\max\left\{\|L\|, \, |\lambda_0|\|\right\}.
\end{align}
Let
\begin{equation*}\alpha=\left\|\prod_{i=1}^n X_i^{(\pi(P),\pi(Q),k_i)}\right\|.
\end{equation*}
By \eqref{norm equation finite product}, \eqref{ley relation for norm-1} and \eqref{ley relation for norm-111},   we have
\begin{align*}\alpha=&\left\|\pi\left(\prod_{i=1}^n X_i^{(P,Q,k_i)}\right)\right\|=\left\|\prod_{i=1}^n \pi\left(X_i^{(P,Q,k_i)}\right)\right\|\\
=&\left\|\prod_{i=1}^n X_i^{(\pi(P),\pi(Q),k_i)}+\widetilde{P_{\mathcal{R}}}\right\|=\max\left\{\alpha, \, \|\widetilde{P_{\mathcal{R}}}\|\right\}=\max\left\{\alpha, 1\right\}.
\end{align*}
Hence $\alpha=1$, which apparently gives  $$\left\|\prod_{i=0}^n X_i^{(\pi(P),\pi(Q),k_i)}\right\|=1$$ by setting
$X_0^{(\pi(P),\pi(Q),k_0)}=I$. It is notable that $\pi(P)$ and $\pi(Q)$ are projections acting on a Hilbert space, so $\big(\pi(P),\pi(Q)\big)$ is harmonious. It follows from Lemmas~\ref{technical lem2-norm} and \ref{technical lem2-norm++} that
\begin{equation*}\|\underline{T}\|\ge |\lambda|\quad\mbox{and}\quad \|Z\|\ge |\lambda_0|.\end{equation*}
So by \eqref{new estimation of L}  we have $\|L\|\ge \|Z\|\ge |\lambda_0|$, which leads by \eqref{new estimation of T down} to  $\|\underline{T}\|=\|L\|$. Combining this equality with $\|\underline{T}\|\ge |\lambda|$ and \eqref{new estimation of T}, we arrive at
$\|T\|=\|\underline{T}\|=\|L\|$.
\end{proof}

Under the restriction of $\lambda_0=0$ in  \eqref{defn of T and tilde T}, a corollary can be derived immediately as follows.
\begin{corollary}\quad Let $C^*(P,Q,P_{\mathcal{R}})$ and $C^*(P-P_{\mathcal{R}}, Q-P_{\mathcal{R}})$ denote the $C^*$-subalgebras of $\mathcal{L}(H)$ generated by elements in $\{P, Q, P_{\mathcal{R}}\}$ and $\{P-P_{\mathcal{R}}, Q-P_{\mathcal{R}}\}$, respectively. Then
each faithful representation $(\pi,X)$  of $C^*(P,Q,P_{\mathcal{R}})$ can induce a faithful representation $(\widetilde{\pi},X)$ of $C^*(P-P_{\mathcal{R}}, Q-P_{\mathcal{R}})$  such that \eqref{two natural map ranges} is satisfied.
\end{corollary}

In the derivations given as above, we merely consider the meaningfulness of $P_{\mathcal{R}}$. At this moment, we don't know whether there exist projections $P$ and $Q$ such that $\mathcal{R}$ is orthogonally complemented, whereas $\mathcal{N}$ fails to be orthogonally complemented\footnote{\,Alternatively, is it possible to find a matched triple $(P,Q,H)$ such that $(I-Q,I-P,H)$ is not a matched triple?}.
To give a partial answer, we need an auxiliary lemma, whose proof is given for the sake of completeness.

\begin{lemma}\quad\label{lem:auxiliary lemma-01}Let $P,Q\in\mathcal{P}(H)$ be such that the sequence $\{(PQP)^n\}_{n=1}^\infty$ converges to $T\in\mathcal{L}(H)$ in norm-topology. Then $T$ is a projection such that $\mathcal{R}(T)=\mathcal{R}$, where $\mathcal{R}$ is defined by \eqref{defn of R and N}.
\end{lemma}
\begin{proof}[\bf Proof]Clearly, $T$ is a projection such that $\mathcal{R}\subseteq \mathcal{R}(T)$ and $PT=T$. So it needs only to show that $QT=T$, or equivalently, $\left[(I-Q)T\right]^*(I-Q)T=0$, which can  be derived from the equations
$$(PQP)^n (I-Q)(PQP)^n=(PQP)^{2n}-(PQP)^{2n+1},\quad\forall n\in\mathbb{N}.$$
\end{proof}

\begin{corollary}\quad \label{cor:be applied} Let $(P,Q,H)$ be a matched triple  such that $\|(P-P_{\mathcal{R}})(Q-P_{\mathcal{R}})\|<1$. Then  $(I-Q,I-P,H)$ is also a matched triple.
\end{corollary}
\begin{proof}[\bf Proof]  Choose any faithful unital representation $(\pi,X)$ of $\mathcal{L}(H)$.   Let $P_{\mathcal{R}\pi}$ be  defined by \eqref{p pi R}, and put
\begin{align}\label{defn of p pi N}&P_{\mathcal{N}\pi}=\big(\pi(I-Q)\big)\wedge \big(\pi(I-P)\big),\\
&S=\left[I-\pi(Q)-P_{\mathcal{N}\pi}\right]\left[I-\pi(P)-P_{\mathcal{N}\pi}\right],\nonumber\\
&T=\left[\pi(P)-P_{\mathcal{R}\pi}\right]\left[\pi(Q)-P_{\mathcal{R}\pi}\right],\nonumber\\
&W=(I-P)(I-Q)(I-P).\nonumber
\end{align}
Then
\begin{align*}P_{\mathcal{N}\pi}=&P_{\mathcal{R}\big(\pi(I-Q)\big)\cap \mathcal{R}\big(\pi(I-P)\big)}=P_{\mathcal{R}\big(I-\pi(Q)\big)\cap \mathcal{R}\big(I-\pi(P)\big)}\\
=&P_{\mathcal{N}(\pi(P))\cap \mathcal{N}(\pi(Q))}.
\end{align*}
By Theorems~\ref{technical thm} and \ref{thm:algebras isomorphic-01}, there exists a unitary $U\in\mathbb{B}(X)$ such that
\begin{align*}\|S\|=&\|U^*TU\|=\|T\|=\|(P-P_{\mathcal{R}})(Q-P_{\mathcal{R}})\|<1,
\end{align*}
hence $\|(S^*S)^n\|=\|S^*S\|^n=\|S\|^{2n}\to 0$ as $n\to\infty$. For each $n\in\mathbb{N}$, it is clear that
$$(S^*S)^n=\big(\pi(W)\big)^n-P_{\mathcal{N}\pi},$$  so $\left\{\big(\pi(W)\big)^n\right\}_{n=1}^\infty$ is norm-convergent and thus is a  Cauchy sequence, and so does for $\{W^n\}_{n=1}^\infty$ by the faithfulness of $\pi$. Due to the completeness of $\mathcal{L}(H)$, $\{W^n\}_{n=1}^\infty$  is norm-convergent.  The assertion then follows from Lemma~\ref{lem:auxiliary lemma-01}.
\end{proof}

\begin{remark}\quad Our next example shows that there exists a matched triple $(P,Q,H)$ such that the sequence $\{(PQP)^n\}_{n=1}^\infty$ does not converge strongly to $P_{\mathcal{R}}$.
\end{remark}
\begin{example}\quad\label{ex:not sot convergence}Let $\mathfrak{A}, H, P$ and $Q$ be as in Example~\ref{ex:counter example-01}. According to \eqref{equ:intersections of null spaces are zeros}, we have $P_{\mathcal{R}}=0$.
From \eqref{equ:defn of P t and Q t} we obtain
\begin{align*}\big[(PQP)^n\big](0)=\left(
                                        \begin{array}{cc}
                                          1 &  \\
                                           & 0 \\
                                        \end{array}
                                      \right)
\end{align*}
for every $n\ge 1$. It follows that
$$\left\|(PQP)^n P\right\|\ge \left\|\left[(PQP)^n P\right](0)\right\|=1, \quad \forall n\ge 1.$$
Thus, $\{(PQP)^n\}_{n=1}^\infty$ does not converge strongly to zero.
\end{example}

\bigskip


\end{document}